\documentclass[12pt,amsfonts, norsk]{article}
\usepackage{graphicx}
\usepackage{latexsym}
\usepackage{amssymb}
\usepackage{amsmath}
\usepackage{layout}
\newtheorem{prop}{Proposition}
\newtheorem{lemma}{Lemma}
\newtheorem{definition}{Definition}
\newtheorem{corollary}{Corollary}
\newtheorem{thm}{Theorem}

\def\real{{\mathord{{\rm I\kern-2.8pt R}}}}        % Fake blackboard bold R.
\def\inte{{\mathord{{\rm I\kern-2.8pt N}}}}

\def\sZZ{{\rm Z\kern-2.8ptem{}Z}}

\def\z{{\mathchoice
  {\sZZ}
  {\sZZ}
  {\rm Z\kern-0.30em{}Z}
  {\rm Z\kern-0.25em{}Z} }}
\def\sQQ{{\kern 0.27em \vrule height1.45ex width0.03em depth0em
          \kern-0.30em \rm Q}}
\def\qu{{\mathchoice
    {\sQQ}
    {\sQQ}
  {\kern 0.225em \vrule height1.05ex width0.025em depth0em \kern-0.25em \rm Q}
  {\kern 0.180em \vrule height0.78ex width0.020em depth0em \kern-0.20em \rm Q}
        }}
\def\sCC{{\kern 0.27em \vrule height1.45ex width0.03em depth0em
          \kern-0.30em \rm C}}
\def\complex{{\mathchoice
    {\sCC}
    {\sCC}
  {\kern 0.225em \vrule height1.05ex width0.025em depth0em \kern-0.25em \rm C}
  {\kern 0.180em \vrule height0.78ex width0.020em depth0em \kern-0.20em \rm C}
        }}

\newcommand{\fin}
{ \vspace{-0.6cm}
\begin{flushright}
\mbox{$\Box$}
\end{flushright}
\noindent }

\newcommand{\R}{\mathbb{R}}

\newcommand{\ba}{\begin{array}}
\newcommand{\ea}{\end{array}}
\newcommand{\be}{\begin{equation}}
\newcommand{\ee}{\end{equation}}
\newcommand{\bea}{\begin{eqnarray}}
\newcommand{\eea}{\end{eqnarray}}
\newcommand{\beaa}{\begin{eqnarray*}}
\newcommand{\eeaa}{\end{eqnarray*}}

\def\z{\zeta}

\font\tenmath=msbm10 \font\sevenmath=msbm7 \font\fivemath=msbm5
\newfam\mathfam \textfont\mathfam=\tenmath
\scriptfont\mathfam=\sevenmath \scriptscriptfont\mathfam=\fivemath

\def \={{\buildrel {\rm (law)} \over =}}

\def \R{{\math R}}

\def\R{{\mathbb{R}}}

\newcommand{\basa}{\begin{assumption}}
\newcommand{\easa}{\end{assumption}}

\newcommand{\bas}{\begin{assum}}
\newcommand{\eas}{\end{assum}}

\newcommand{\ignore}[1]{}
\begin{document}

\renewcommand{\thefootnote}{\fnsymbol{footnote}}

\title{{\bf Optimal control for rough differential equations}}
\date{
\normalsize
{\bf Laurent Mazliak\footnote{Corresponding author. Fax: +33 1 47 27 72 23}\,\, and Ivan Nourdin}\\
\small{
Universit\'e Pierre et Marie Curie (Paris 6)\\
Laboratoire Probabilit\'es et Mod\`eles Al\'eatoires (LPMA)\\
Bo\^ite courrier 188, F-75252 Paris Cedex 5\\
{\tt \{mazliak,nourdin\}@ccr.jussieu.fr}
}
}
\maketitle

\begin{abstract}

\vspace{2mm}

\noindent
In this note, we consider an optimal control problem associated to a differential 
equation driven by a H\"older continuous function $g$ of index $\beta>1/2$.
We split our study in two cases. If the coefficient of $dg_t$ does not
depend on the control process, we prove an existence theorem for a slightly
generalized control problem, that is we obtain a literal extension of 
the corresponding deterministic situation. If the coefficient of $dg_t$ depends on
the control process, we also prove an existence theorem but we are here obliged to 
restrict the set of controls to sufficiently regular functions.\\
\end{abstract}


\section{Introduction}

In recent years, several authors have been interested in control problems involving a stochastic process driven
by a fractional Brownian motion (fBm in short). This kind of situation may indeed be natural when one wishes to
modelized a problem in which  long distance memory effects may occur. However, contrary to the situation of
processes driven by ordinary Brownian motion, fBm lacks many strong stochastic properties
(in particular, {\it martingales} properties). Although attempts have been made to build 
a stochastic calculus for fBm some twenty years ago (see \cite{Coutin} for a nice survey), 
the developed techniques remain until today not so easy to use.

In order to deal with the situation of control for processes driven with fBm, it was up to now necessary to limit oneself
to special situations. Two directions have been mainly considered in recent literature:
\begin{enumerate}
\item A series of paper used a
stochastic calculus based on Wick integrals and developed in \cite{Duncan}; see in particular \cite{HuOksendalSulem}
for considerations about application to mathematical finance. Moreover, in a recent paper \cite{HuZhou}, the special
situation of a linear regulator driven by fBm was considered and optimal controls characterized.
\item A completely different way has been followed in several papers of a group in France, but only when the considered problem
is strongly linear, a situation permitting to the author to use a martingale representation of fBm in order to transform
the original problem in an ordinary stochastic control problem (see in particular \cite{Kleptsyna} for a general exposition
of the method).
\end{enumerate}

In the present paper, we also keep a stochastic control problem for a process driven by a fBm in  the background, and we want to study existence results for optimal controls. Nevertheless, as the basic property we need in our approach occurs to be the H\"older regularity of the paths of fBm, we study here the deterministic situation
where the state under control is driven by a H\"older continuous function $g$ of index $\beta$
\begin{equation}\label{edsintro}
x^u_t=x^u_0+\int_0^t \sigma(s,u_s,x^u_s)dg_s + \int_0^t b(s,u_s,x^u_s)ds,\quad t\in[0,T],
\end{equation}
where the control process $u$ belongs to a set of admissible controls ${\cal U}$.
When $\beta\in(1/2,1)$, it is possible to choose Young integral \cite{young} for integration with respect to $dg_t$ in (\ref{edsintro}), which simply appears to be the limit of Riemann sums.
Moreover, as it was remarked in \cite{zahle}, it is possible to express it in terms of fractional derivative operators
(see Section 2 below). This appears to be very useful to allow easier computations, and mostly 
 has motivated us to restrict to the case where $\beta\in(1/2,1)$ in order to choose the Young integral. The next step would be considering the case where $\beta\in(1/3,1/2]$. In the very
recent work \cite{HN2}, an equivalent to the expression of the Young integral
in terms of fractional derivative operators has been proposed, using quadratic multiplicative functional. Thus,
the strategy introduced in this paper could also certainly be derived for $\beta\in(1/3,1/2]$ and
we propose to study this fact in a forthcoming paper. 

The control problem considered in the present paper can be set in the following way. 
\noindent  \begin{equation}\label{question}
\begin{array}{lll}
\mbox{{\bf Problem}: `A cost functional $J:{\cal U}\rightarrow\R$ being given,
is it possible to prove}\\
\mbox{the existence of $u^*\in{\cal U}$ realizing $\inf_{u\in{\cal U}} J(u)$?'}
\end{array}
\end{equation}

As usual, the bigger ${\cal U}$ is, the more difficult it is to answer this question.  A general methodology is to look for conditions ensuring that ${\cal U}$ be compact for a certain topology under which $J$  is continuous.

Differential equations of the type (\ref{edsintro}) (without the control process $u$) have been intensively studied in recent years, due in particular to possible applications for fBm (see, {\it e.g.}, \cite{HN,HN2,lyons,NS2,ruz}).
Since the obtention of solutions to (\ref{edsintro}) requires in general regularity on the coefficients
(see Theorem \ref{thm1} below), we split our study in two cases. 
\begin{enumerate}
\item If the coefficient of $dg_t$ does {\it not} depend on the control process,
we are able to extend the deterministic
situation, and to prove an existence theorem for a slightly generalized control problem where the controls are in
fact randomized : see Corollary \ref{cor1} and Proposition \ref{prop4}. In fact, we use the so-called `compactification methods', which have been developed during the 1960's
for deterministic control problems (see  \cite{Ghoula}, \cite{Young1969}), and during the 1970's for the stochastic control problem
(see \cite{Fleming}, \cite{ElKaroui}).
\item  If the coefficient of $dg_t$ {\it does } depend  on the control function,
the situation is much more intricate, and
this obliges us to severely restrict the set of controls to sufficiently regular functions.
A challenging question would
be to relax this hypothesis, but this would require to get a reasonable notion of solution for a differential equation
with a weaker regularity than H\"olderian. This seems not to be already available in the literature, up to our best knowledge.
\end{enumerate}

The paper is organized as follows. In section 2, we recall some now classical facts on fractional integrals
and derivatives, which are useful for the sequel. In section 3, we study the optimal control problem in the case where $\sigma$ does not depend on $u$. The case where $\sigma$ depends on $u$ is considered in section 4.

\section{Fractional integrals and derivatives}

Let $a,b\in\R$, $a<b$. For any $p\geq 1$, we denote by $L^p=L^p([a,b])$
the usual Lebesgue spaces of functions on
$[a,b]$.

Let $f\in L^1$ and $a>0$. The left-sided and right-sided fractional Riemann-Liouville
integrals of $f$ of order $\alpha$ are defined for almost all $x\in
(a,b)$ by
$$
 I^{\alpha}_{a+}f(x)  = \frac{1}{\Gamma(\alpha)}\int_{a}^x
 (x-y)^{\alpha-1}f(y)dy,
$$
and
$$
 I^{\alpha}_{b-}f(x)  = \frac{(-1)^{-\alpha}}{\Gamma(\alpha)}\int_{x}^b
 (y-x)^{\alpha-1}f(y)dy,
$$
respectively,
where $(-1)^{-\alpha}={\rm e}^{-i\pi\alpha}$ and $\Gamma(\alpha)=\int_0^\infty r^{\alpha-1}{\rm e}^{-r}dr$
denotes the usual Euler function.

If $f\in I^{\alpha}_{a+}(L^p)$ (resp. $f\in I^{\alpha}_{b-}(L^p)$)
and $\alpha\in
(0,1)$, then for almost all $x\in (a,b)$, the left-sided and right-sided
Riemann-Liouville derivative of $f$ of order $\alpha$ are defined by
\begin{equation}\label{d+}
D^{\alpha}_{a+}f(x) = \frac{1}{\Gamma(1-\alpha)}\left ( \frac{f(x)}{(x-a)^{\alpha}}+
        \alpha
\int_{a}^x \frac{f(x)-f(y)}{(x-y)^{\alpha+1}}dy \right)
\end{equation}
and
\begin{equation}\label{d-}
D^{\alpha}_{b-}f(x) = \frac{1}{\Gamma(1-\alpha)}\left ( \frac{f(x)}{(b-x)^{\alpha}}+
        \alpha
\int_{x}^b \frac{f(x)-f(y)}{(y-x)^{\alpha+1}}dy \right)
\end{equation}
respectively, where $a\le x\le b$.

If $\mu\in(0,1)$ and $T\in (0,\infty)$, we note ${\cal C}^\mu([0,T])$ the set of functions
$g:[0,T]\rightarrow\R$ such that
$$
\sup_{0\le s< t\le T}\frac{|g(t)-g(s)|}{|t-s|^\mu}<+\infty.
$$
If there is no ambiguity, we prefer note ${\cal C}^\mu$
instead of ${\cal C}^\mu([0,T])$.
The set ${\cal C}^\mu$ is a Banach space when it is endowed with the following norm:
$$
|g|_{\infty,\mu} := \sup_{0\le t\le T}|g(t)| + \sup_{0\le s< t\le T}\frac{|g(t)-g(s)|}{|t-s|^\mu}.
$$
We also set, for $a,b\in [0,T]$ and $g\in{\rm C}^\mu$:
$$
|g|_{a,b,\mu}= \sup_{a\le s< t\le b}\frac{|g(t)-g(s)|}{|t-s|^\mu}.
$$
and
$$
|g|_{a,b,\infty}=\sup_{a\le t\le b}|g(t)|.
$$
When $a=0$ and $b=T$ we simply note $|g|_\mu$ and $|g|_\infty$ instead of $|g|_{0,T,\mu}$ and
$|g|_{0,T,\infty}$ respectively.

Let $f:\R\rightarrow\R\in {\rm C}^\lambda$ and $g:\R\rightarrow\R\in{\rm C}^\mu$ with $\lambda,\mu\in (0,1)$
such that $\lambda+\mu>1$. Then, for any $a,b\in[0,T]$, the Young integral \cite{young}
$\int_a^b f dg$ exists and we can express it in terms of fractional derivatives (see \cite{zahle}):
for any $\alpha\in(1-\mu,\lambda)$, we have
\begin{equation}\label{zahle}
\int_a^b f dg = (-1)^\alpha \int_a^b D_{a+}^\alpha f(x) D_{b-}^{1-\alpha}g_{b-}(x)dx,
\end{equation}
where $g_{b-}(x)=g(x)-g(b)$.

\section{First case: when $\sigma$ does not depend of $u$}
In the sequel, we fix $x_0\in\R$, $\beta\in (1/2,1)$, $T\in (0,\infty)$ and
$g\in{\cal C}^\beta={\cal C}^\beta([0,T])$. We assume moreover that
$\sigma:[0,T]\times\R\rightarrow\R$ is ${\rm C}^{1,2}$ with bounded derivatives
and that $b:[0,T]\times\R^2\rightarrow\R$ is bounded and global Lipschitz, uniformly
in $x\in\R$ with respect to $(t,u)\in[0,T]\times\R$.
\begin{thm}\label{thm5}
For any measurable control $u:[0,T]\rightarrow\R$, the integral equation
\begin{equation}\label{eds2}
x^u_t=x^u_0+\int_0^t \sigma(r,x^u_r)dg_r+\int_0^t b(r,x^u_r,u_r)dr,\quad t\in [0,T]
\end{equation}
admits a unique solution $x^u\in {\cal C}^0([0,T])$.
\end{thm}
{\bf Proof of Theorem \ref{thm5}}.
\begin{itemize}
\item
We first prove Theorem \ref{thm5} in the autonomous case,
that is when $\sigma(t,x)=\sigma(x)$ and $b(t,x,u)=b(x,u)$.
In other words,
we consider
\begin{equation}\label{eds21}
x^u_t=x^u_0+\int_0^t \sigma(x^u_r)dg_r+\int_0^t b(x^u_r,u_r)dr,\quad t\in [0,T]
\end{equation}
instead of (\ref{eds2}).
At this level, we need a preliminary lemma:
\begin{lemma}\label{lm8}
Assume that $h:[0,T]\times\R^3\rightarrow\R$ is such that, for any $R>0$, there exists
$c_R>0$ verifying
\begin{equation}\label{f}
\forall (r,g,u,y,z)\in [0,T]\times[-R,R]\times \R^3:\quad |h(r,g,u,y)-h(r,g,u,z)|\le
c_R|y-z|,
\end{equation}
and that $u:[0,T]\rightarrow\R$ is a measurable function.
Then the integral equation
\begin{equation}\label{eds14}
y_t=y_0+\int_0^t h(r,g_r,u_r,y_r)dr,\quad t\in [0,T]
\end{equation}
admits a unique solution $y\in{\rm C}^{0}([0,T])$.
\end{lemma}
{\bf Proof of Lemma \ref{lm8}}. We only sketch the proof, the arguments used being
classical.\\
{\it Existence}. Let us define $(y^n)$ recursively by $y^0(t)\equiv y_0$ and
$$
y^{n+1}(t)=y_0+\int_0^t h(r,g(r),u(r),y^n(r))dr,\quad t\in [0,T].
$$
Since $g$ is continuous, there exists $R>0$ such that $g([0,T])\subset [-R,R]$.
Thus, using the hypothesis made on $h$, it is classical to prove that
$|y^{n+1}-y^n|_\infty\le \frac{{c_R}^n}{n!}$. In particular, the sequence
$(y^n)$ is Cauchy and the limit $y$ is a solution to (\ref{eds14}).\\
{\it Uniqueness}. Let $y$ and $z$ be two solutions of (\ref{eds14}). Then,
for any $t\in [0,T]$, we easily have
$$
|y-z|_{\infty,[0,t]}\le c_R\int_0^t |y-z|_{\infty,[0,r]}dr
$$
and we can conclude that $y=z$ using Gronwall's lemma.
\fin
We now apply the Doss-Sussmann's method in order to finish the proof
of Theorem \ref{thm5} in the autonomous case. First, we denote by $\phi$ the unique
solution to
\begin{equation}\label{phi}
\frac{\partial\phi}{\partial g}(g,y)=\sigma\circ\phi(g,y),\,\forall g,y\in\R
\mbox{ and }\phi(0,y)=y,\,\forall y\in\R.
\end{equation}
The hypothesis made on $\sigma$ ensures that $\phi$ is well-defined.
We also have, for $g,y\in\R$:
$$\frac{\partial\phi}{\partial y}(g,y)={\rm exp}\left(
\int_0^g \sigma'(\phi(h,y)) dh
\right).$$
Define $f:\R^3\rightarrow\R$ by
$$f(g,u,y)=\frac{b(\phi(g,y),u)}{\frac{\partial\phi}{\partial y}(g,y)}
=b(\phi(g,y),u)\,{\rm exp}\left(-
\int_0^g \sigma'(\phi(\ell,y)) d\ell\right).
$$
The
hypothesis made on $b$ and $\sigma$ ensures that $h:[0,T]\times\R^3\rightarrow\R$
defined by $h(r,g,u,y)=f(g,u,y)$
verifies (\ref{f}). Thus, there exist a unique $y$ solution to (\ref{eds14}).
Using the change of variable formula, it is now immediate to prove
that $x^u_t=\phi(g_t,y_t)$ is a solution to (\ref{eds21}).
For the uniqueness, it suffices to adapt to our context the
proof contained in \cite{Doss}, page 103.
\item The general case being similar with the previous case, we only sketch the proof. Here, we have
to consider $\phi$ given by
\begin{equation}\label{phigen}
\frac{\partial\phi}{\partial g}(r,g,y)=\sigma(r,\phi(r,g,y)),\,\forall (r,g,y)\in[0,T]\times\R^2
\end{equation}
with initial conditions
$$
\phi(r,0,y)=y,\,\forall (r,y)\in[0,T]\times\R
$$
instead of (\ref{phi}). Moreover, $y:[0,T]\rightarrow\R$ is, in the case,
defined as the unique solution to (\ref{eds14}) with $h$ given by
$$
h(r,g,u,y)=\frac{b(r,\phi(r,g,y),u)-\frac{\partial\phi}{\partial r}(r,g,y)}
{\frac{\partial\phi}{\partial y}(r,g,y)},
$$
see also \cite{Doss}, page 116.
Finally, the unique solution to (\ref{eds2}) is given by
$$x^u_t=\phi(t,g_t,y_t).$$
\end{itemize}
\fin

In order to make use of a compactification method, it is necessary to enlarge the set of controls by considering {\it relaxed controls}.

\begin{definition}
A relaxed control is a measure $q$ over $U\times [0,T]$ such that the projection of $q$
on $[0,T]$ is the Lebesgue measure.
We denote by $\cal V$ the set of relaxed controls.
\end{definition}

A relaxed control $q$ can be decomposed with a measurable kernel: $q(da,dt)=q_t(da)dt$ where $t\mapsto q_t$ is a measurable
function from $\R^+$ to the set of probability measures on $U$. There is a natural embedding of (non-relaxed) controls in
the set of relaxed controls: $q$ is a non-relaxed control if at each time $t$, $q_t$ concentrates on a single point $u_t$.
In other words, we assimilate the control $(u_t)_{t\in [0,T]}$ with the relaxed control $\delta_{u_t}dt$ where $\delta_x$
denotes the Dirac measure at $x$. We denote by $\cal V^0$ the set of non-relaxed controls.

The main result that we shall need is the immediate following consequence of the vague topology.

\begin{prop} Suppose $U$ is a compact subset of $\R^n$. The set $\cal V$ of relaxed controls equipped with the vague
topology is compact.
\end{prop}

From now on, we shall suppose that the set $U$ is compact. A solution to equation (\ref{eds2}) associated to a relaxed control $q$ is obtained in the following extension of Theorem \ref{thm5}.

\begin{thm}\label{thm6} Let $q\in \cal V$ be a relaxed control. There exists a unique solution $x^q \in C^0([0,T])$ of the equation
\begin{equation}\label{star}
x_t^q=x_0^q+\int_0^t\sigma (r,x_r^q)dg_r+\int_0^t\int_U b(r,x_r^q,a)q_r(da)dr.
\end{equation}
Moreover, $q\mapsto x^q$ is continuous from $\cal V$ to $C^0([0,T])$.
\end{thm}
{\bf Proof}.
Denote by $\phi$ the unique solution to (\ref{phigen}).
Set
\begin{equation}\label{h}
h(r,g,q,y)={\int_U b(r,\phi (r,g,y),a)q_r(da)- {\partial \phi \over \partial r}(r,g,y) \over
{\partial \phi \over \partial y}(r,g,y)}.
\end{equation}
Clearly, due to the hypotheses on $b$ and $\sigma$, $\forall (r,g,q,y,z)\in
[0,T]\times[-R,R]\times {\cal V}\times \R \times \R$,
$$|h(r,g,q,y)-h(r,g,q,z)|\le c_R |y-z|.$$
Therefore, the integral equation
(\ref{eds14}) admits a
unique solution $y\in C^0([0,T])$, see Lemma \ref{lm8}.
Then, one may check that $x_t^q =\phi (t,g_t,y_t)$ is a solution to (\ref{star}).
Uniqueness is obtained as before.

Suppose now that $q^n$ is a sequence in $\cal V$, converging to $q\in \cal V$
%Set $$h_n(r,g,q^n,y)={\int_U b(\phi (g,y),a)q^n_r(da)
%\over {\partial \phi \over \partial y}
%(g,y)}
%=\int_U b(\phi (g,y),a)q^n_r(da)\,{\rm exp}\left(-
%\int_0^g \sigma'(\phi(\ell,y)) d\ell\right).
%$$
and let $y^n$ be the solution of (\ref{eds14}) associated to $h=h(r,g,q^n,y)$ given by (\ref{h}).
Using hypotheses on $b$, we now prove that $y^n$
converges to $y$ in $C^0([0,T])$. Indeed,
$$
\begin{array}{lll}
|y_t-y^n_t|&=&\left|\int_0^t[h(s,g_s,q,y_s)-h(s,g_s,q^n,y^n_s)]ds\right|\\
&\le& \left|\int_0^th(s,g_s,q,y_s)ds-\int_0^t h(s,g_s,q^n,y_s)ds\right|\\
&&+\left|\int_0^t h(s,g_s,q^n,y_s)ds-\int_0^t h(s,g_s,q^n,y^n_s)ds\right|\\
&\le&\left|\int_0^t \int_U b(s,\phi (s,g_s,y_s),a)q^n_s(da)ds- \int_0^t \int_U  b(s,\phi (s,g_s,y_s),a)q_s(da)ds\right|\\
&&+c_R \int_0^t |y_s-y^n_s|ds .
\end{array}
$$

In the last expression, the first term tends to 0 due to the vague convergence of $q^n_s(da)ds$ to $q_s(da)ds$,
and continuity and boundedness hypotheses on $b$.
It results therefore from Gronwall's lemma that $|y-y^n|_{\infty}$ tends to 0.
Finally, as the solution $x^q$ (resp. $x^{q^n}$) of (\ref{star}) associated to $q$ (resp. $q^n$)
is given by $x_t =\phi (t,g_t,y_t)$ (resp. $x_t^n =\phi (t,g_t,y^n_t)$), one easily deduces that $|x-x^n|_{\infty}$
tends to 0.
\fin

Consider now a cost in  integral form: for $u_t$ a given control taking values in $U$, we set
$$J(u)=\int_0^T \ell (r,x_r^u,u_r)dr$$ where $\ell$ is a bounded continuous function on
$(r,x,u)$. The definition can be  immediately extended to the case of relaxed controls: if $q$ is a relaxed control from $\cal V$, $$J(q)=\int_0^T \int_U\ell (r,x_r^q,a)q_r(da)dr.$$ Using the continuity property of Theorem \ref{thm6}, and the hypotheses on $\ell$, one obtains the following Proposition.

\begin{prop} Under the hypotheses introduced in the previous paragraph, the
application $q\mapsto J(q)$ is continuous on $\cal V$. \end{prop}

The set  $\cal V$ being compact, one immediately deduces the following existence
result.

\begin{corollary}\label{cor1} Under the prevailing hypotheses, there exists $q^*\in \cal V$ such
that $$J(q^*)=\inf_{q\in\cal V}J(q).$$
\end{corollary}

We conclude the present section by proving that one has not enlarged too much the control problem by considering relaxed controls. More precisely, we now prove that the optimal cost (i.e. the infemum of the cost functional) over relaxed and non-relaxed controls is the same.  This result is obtained as in the deterministic case by means of approximation of relaxed controls by {\it step constant} relaxed controls, and then by non-relaxed controls via the so-called {\it chattering lemma}, a method originally introduced in \cite{Ghoula}. We here only sketch these two steps. \vskip 3mm
{\it First step} : $q\in \cal V$ is approximated by relaxed controls of the form
$$\sum_{j=0}^{N-1}\sum_{i=1}^kq_i^j\delta_{a_i}(da){\bf 1}_{[t_j,t_{j+1}[}$$ where
$0=t_0<t_1<\dots <t_N=T$, $a_1,\dots ,a_k$ are elements in $U$, and for each $j=0,\dots ,N-1$, $q_1^j,\dots ,q_k^j$ are non-negative real numbers such that $\displaystyle \sum_{i=1}^kq_i^j=1$. This is a straightforward consequence of approximation of the measurable function $t\mapsto q_t$ by a step function and of approximation of a probability measure $\mu$ on $U$ by point measures of the form
$\displaystyle \sum_{i=1}^m \mu_i \delta_{a_i}$.
\vskip 2mm
{\it Second step}:  Recall the chattering lemma (see \cite{Ghoula}, Theorem 1)

\begin{prop} Let $a_1,\dots ,a_k$ be in $U$ and $q_1,\dots ,q_k$ be non-negative real numbers such that $\displaystyle \sum_{i=1}^kq_i=1$. Let $f$ be a bounded continuous function from $[s,t]\times U$ to $\R$. Then, for $\varepsilon >0$ given, there exists a measurable partition $V_1,\dots ,V_k$ of $[s,t]$ such that
$$\left|\int_s^t\sum_{i=1}^kq_if(r,a_i)dr-\sum_{i=1}^k\int_{V_i}f(r,a_i)dr\right|<\varepsilon .$$
\end{prop}
In other words, the step-relaxed control $\big(\sum_{i=1}^kq_i\delta_{a_i}(da)\big)dr$ is approximated by the non-relaxed
control $\sum_{i=1}^k{\bf 1}_{V_i}(r)\delta_{a_i}(da)dr.$
Therefore, for any $q\in \cal V$, there exists a sequence of (non-relaxed) controls $q^n$ which converges to $q$.
As, obviously, the infemum of $J$ on $\cal V$ is smaller than the infemum on $\cal V^0$, using the continuity of $J$, we obtain the following comparison result.

\begin{prop}\label{prop4} Under the hypotheses of the present section,  $\inf_{q\in \cal V}J(q)=\inf_{q\in \cal V^0}J(q) .$
\end{prop}

\section{Second case: when $\sigma$ depends of $u$}

As already mentioned in the introduction, the case when $u$ enters the coefficient of $dg_t$ seems to be much more complicated as we do not have a reasonable way for integrating functions less regular than H\"olderian. Therefore we shall need to restrict very strongly our admissible controls set.

In the sequel, we fix $x_0\in\R$, $\beta\in (1/2,1)$, $T\in (0,\infty)$,
$g\in{\cal C}^\beta={\cal C}^\beta([0,T])$,
$\sigma:[0,T]\times\R^2\rightarrow\R\in{\rm C}^{1,2,2}$ with bounded derivatives
and $b:[0,T]\times\R^2\rightarrow\R$ global Lipschitz continuous.
\begin{thm}\label{thm1}
For any control $u\in{\cal C}^\mu$ with $1-\beta<\mu\le \beta$, the integral equation
\begin{equation}\label{eds}
x^u_t=x^u_0+\int_0^t \sigma(r,x^u_r,u_r)dg_r+\int_0^t b(r,x^u_r,u_r)dr,\quad t\in [0,T]
\end{equation}
admits a unique solution $x^u\in {\cal C}^\mu$.
\end{thm}
{\bf Proof of Theorem \ref{thm1}}. It suffices to adapt the proof of Ruzmaikina \cite{ruz} to our context,
{\it i.e.} to hold account of the control $u$. There is not new difficulties. See also Nualart and
R${\check {\rm a}}$s\c{c}anu \cite{NR}.
\fin
\begin{thm}\label{thm4}
If $\cal U$ is a set of functions which is bounded in a certain ${\cal C}^\mu$ with $\mu\in (1-\beta,\beta]$
and if $J:\cal U\rightarrow \R$ is continuous for $|\cdot|_{\infty,\mu'}$ for a certain $\mu'\in(1-\beta,\mu)$
then the following control problem can be solved:
$$\mbox{{\rm there exists $u^*\in{\cal U}$ realizing $\inf_{u\in{\cal U}} J(u)$}}.$$
\end{thm}
{\bf Proof of Theorem \ref{thm4}}. It is a direct consequence of Lemma \ref{thm3} below.\fin
\begin{lemma}\label{thm3}
If $\cal U$ is a set of functions which is bounded in a certain ${\cal C}^\mu$ with $\mu\in (1-\beta,\beta]$
then the set of all couples $(u,x^u)\in {\cal U}\times{\cal C}^\mu$
is relatively compact in ${\cal C}^{\mu'}\times {\cal C}^{\mu'}$ for any $\mu'\in(1-\beta,\mu)$.
\end{lemma}
{\bf Proof of Lemma \ref{thm3}}. According to Lamperti \cite{lamperti}, we know
that $\cal U$ is relatively compact in ${\cal C}^{\mu'}$.
Lemma \ref{thm2} below allows then to conclude.\fin
\begin{lemma}\label{thm2}
For any $\mu\in (1-\beta,\beta]$, $T^\mu:{\cal C}^\mu\rightarrow{\cal C}^\mu$ defined by
$T^\mu(u)=x^u$ is a continuous operator.
\end{lemma}
{\bf Proof of Lemma \ref{thm2}}. We adapt the proof of Theorem 3.2 in Hu and Nualart \cite{HN}.
For simplicity, we assume that $\sigma(r,x,u)=\sigma(x,u)$
and $b(r,x,u)=b(x,u)$, the proof of the general case being similar.
Moreover positive constants, depending only on $b$, $\sigma$, their derivatives, $x_0$ and $g$, will be
denoted by $k$, regardless of their value.
Let $(u^n)\subset {\cal C}^\beta$ be such that $u^n\rightarrow u$.
%We denote $x^{u_n}=T^\mu(u^n)$ and $x^u=T^\mu(u)$.
Fix $s,t\in[0,T]$ and let $\alpha\in (1-\beta,\mu)$. We can write, using (\ref{zahle}):
$$
\begin{array}{lll}
&|x^u_t-x^{u_n}_t-x^u_s+x^{u_n}_s|\\
=&\left|\int_s^t [\sigma(x^u_r,u_r)-\sigma(x^{u_n}_r,u^n_r)]dg_r+
\int_s^t [b(x^u_r,u_r)-b(x^{u_n}_r,u^n_r)]dr\right|\\
\le&\int_s^t |D_{s+}^\alpha[\sigma(x^u_r,u_r)-\sigma(x^{u_n}_r,u^n_r)]|\cdot|D_{t-}^{1-\alpha}g_{t-}(r)|dr
+\int_s^t |b(x^u_r,u_r)-b(x^{u_n}_r,u^n_r)|dr.
\end{array}
$$
Using (\ref{d-}), it is easy, on one hand, to show that
$$
|D_{t-}^{1-\alpha} g_{t-}(r)|\le k |g|_\beta |t-r|^{\alpha+\beta-1}.
$$
On the other hand, we have, using (\ref{d+}):
$$
\begin{array}{lll}
&|D_{s+}^\alpha[\sigma(x^u_r,u_r)-\sigma(x^{u_n}_r,u^n_r)]|\\
\le&
|\sigma'|_\infty \left(
|x^u-x^{u_n}|_{s,t,\infty}+|u-u^n|_{s,t,\infty}
\right)(r-s)^{-\alpha}\\
+&|\sigma'|_\infty\left(
|x^u-x^{u_n}|_{s,t,\mu}+|u-u^n|_{s,t,\mu}
\right) (r-s)^{\mu-\alpha}\\
+&|\sigma''|_\infty \left(
|x^u-x^{u_n}|_{s,t,\infty}+|u-u^n|_{s,t,\infty}
\right)\left(
|x^{u_n}|_{s,t,\mu}+|u^n|_{s,t,\mu}
\right)
(r-s)^{\mu-\alpha}.
\end{array}
$$
We deduce that
$$
\begin{array}{lll}
&|x^u-x^{u_n}|_{s,t,\mu}\\
\le&
k\left[
(|x^u-x^{u_n}|_{s,t,\infty}+|u-u^n|_{s,t,\infty})(t-s)^{\beta-\mu}
+(|x^u-x^{u_n}|_{s,t,\mu}+|u-u^n|_{s,t,\mu})(t-s)^\beta\right.\\
+&\left.(|x^u-x^{u_n}|_{s,t,\infty}+|u-u^n|_{s,t,\infty})(|x^{u_n}|_{s,t,\mu}+|u^n|_{s,t,\mu})(t-s)^\beta
\right]
\end{array}
$$
and, by rearranging:
\begin{equation}\label{rea}
\begin{array}{lll}
&|x^u-x^{u_n}|_{s,t,\mu}\\
\le&
k(1-k(t-s)^\beta)^{-1}\left[
|x^u-x^{u_n}|_{s,t,\infty}+|u-u^n|_{s,t,\infty}
+|u-u^n|_{s,t,\mu}(t-s)^\beta\right.\\
+&\left.(|x^u-x^{u_n}|_{s,t,\infty}+|u-u^n|_{s,t,\infty})(|x^{u_n}|_{s,t,\mu}+|u^n|_{s,t,\mu})(t-s)^\beta
\right].
\end{array}
\end{equation}
Then, if we set $\Delta=t-s$:
$$
\begin{array}{lll}
&|x^u-x^{u_n}|_{s,t,\infty}\\
\le &|x^u_s- x^{u_n}_s|+|x^u-x^{u_n}|_{s,t,\mu}\Delta^\beta\\
\le&|x^u_s- x^{u_n}_s|+
k(1-k\Delta^\beta)^{-1}\Delta^\beta\left[
|x^u-x^{u_n}|_{s,t,\infty}+|u-u^n|_{s,t,\infty}\right.\\
+&|u-u^n|_{s,t,\mu}\Delta^\beta
+\left.(|x^u-x^{u_n}|_{s,t,\infty}+|u-u^n|_{s,t,\infty})(
|x^{u_n}|_{s,t,\mu}+|u^n|_{s,t,\mu}
)
\Delta^\beta
\right].
\end{array}
$$
By rearranging, we obtain
$$
\begin{array}{lll}
|x^u-x^{u_n}|_{s,t,\infty}
\le\left[1-k(1-k\Delta^\beta)^{-1}\Delta^\beta(1+(|x^{u_n}|_{s,t,\mu}+|u^n|_{s,t,\mu})\Delta^\beta)
\right]^{-1}\\
\times\left(|x^u_s- x^{u_n}_s|+
k(1-k\Delta^\beta)^{-1}\Delta^\beta\left[
|u-u^n|_{\infty,\mu}
+|u-u^n|_{\infty}(|x^{u_n}|_{\mu}+|u^n|_{\mu})
\right]\right).
\end{array}
$$
We can finish as in \cite{HN} to obtain that $|x^u-x^{u_n}|_\infty\rightarrow 0$ as
$n\rightarrow \infty$. Using finally (\ref{rea}), we obtain that
$|x^u-x^{u_n}|_{\infty,\mu}\rightarrow 0$ as
$n\rightarrow \infty$, that is $x^{u_n}\rightarrow x^u$ in ${\cal C}^\mu$. In other words,
$T^\mu$ is a continuous operator from ${\cal C}^\mu$ to himself.
\fin

An example of a cost $J$ satisfying conditions of Theorem \ref{thm4} is
$$J(u)=\int_0^T \ell(r,x_r^u,u_r)dr$$
with $\ell:[0,T]\times\R^2\rightarrow\R$
verifying
$$
\forall (r,x,y,u,v)\in [0,T]\times\R^4,\quad
|\ell(r,x,u)-\ell(r,y,v)|\le {\rm cst}\left(|x-y|+|u-v|\right).
$$

\end{document}